\documentclass[12pt]{amsart}
\usepackage{amsxtra,latexsym,amssymb}
\usepackage{epsfig,rotate,amsthm}
\usepackage{hyperref}

\def\qed{{\hfill $\Box$}}

\def\Z{{\mathbb Z}}

\def\K{{\mathbb K}}

\def\A{\mathcal{A}_{p}(\lambda, \mu, \K_{q}[s, t])}
\def\D{\mathcal{A}_{p_{1}}(\lambda_{1}, \mu_{1}, \K_{q_{1}}[s, t])}
\def \E{\mathcal{A}_{p_{2}}(\lambda_{2}, \mu_{2}, \K_{q_{2}}[s, t])}
\theoremstyle{theorem}
\newtheorem{thm}{Theorem}[section]
\newtheorem{cor}{Corollary}[section]
\newtheorem{prop}{Proposition}[section]

\newtheorem{lem}{Lemma}[section]
\theoremstyle{definition}
\newtheorem{defn}{Definition}[section]
\theoremstyle{remark}

\newtheorem{rem}{\bf Remark}[section]

\begin{document}
\title[Automorphisms for Generalized Weyl Algebras]{The Automorphism Groups for a family of Generalized Weyl Algebras}
\author[X. Tang]{ Xin Tang}
\address{
Department of Mathematics \& Computer Science\\
Fayetteville State University\\
1200 Murchison Road, Fayetteville, NC 28301}
\email{xtang@uncfsu.edu} 
\thanks{This project is supported by an HBCU Master's Degree STEM Program Research Mini-Grant funded by Title III at Fayetteville State University.}
\date{\today}
\begin{abstract}
In this paper, we study a family of generalized Weyl algebras $\{\A\}$ and their polynomial extensions. We will show that the algebra $\A$ has a simple localization $\A_{\mathbb{S}}$ when none of $p$ and $q$ is a root of unity. As an application, we determine all the height-one prime ideals and the center for $\A$, and prove that $\A$ is cancellative. Then we will determine the automorphism group and solve the isomorphism problem for the generalized Weyl algebras $\A$ and their polynomial extensions in the case where none of $p$ and $q$ is a root of unity. We will establish a quantum analogue of the Dixmier conjecture and compute the automorphism group for the simple localization $(\mathcal{A}_{p}(1, 1, \K_{q}[s, t]))_{\mathbb{S}}$. Moreover, we will completely determine the automorphism group for the algebra $\mathcal{A}_{p}(1, 1, \K_{q}[s, t])$ and its polynomial extension when $p\neq 1$ and $q\neq 1$.
\end{abstract}
\maketitle 
\section*{Introduction}
Generalized Weyl algebras (denoted by $R(\sigma, a)$ and called GWAs for short) have been extensively studied in the literature since they were first introduced by Bavula in \cite{Bavula}. Note that many important examples of generalized Weyl algebras have appeared as families of algebras defined via parameters. Thus it is an important question to study the isomorphism classification for a given family of generalized Weyl algebras. One practical way to obtain such an isomorphism classification is to first determine the automorphism group for each algebra in the family. However, the complete determination of the automorphism group for a given algebra has been proved to be extremely difficult, despite the recent significant progress made in the works \cite{CPWZ1, CPWZ2, GY, SU, Yakimov1, Yakimov2}. 

The study of the algebra automorphism group and isomorphism problem for generalized Weyl algebras was initiated by Bavula and Jordan in \cite{BJ}, where they successfully determined the automorphism group and solved the isomorphism problem for several families of rank-one generalized Weyl algebras over the base ring $R=\K[h]$ or $\K[h^{\pm 1}]$ and related algebras. Motivated by the work of Bavula and Jordan, the automorphism group and isomorphism problem have been further studied for some rank-one generalized Weyl algebras in \cite{RS, SV}. Recently, Goodearl and Hartwig have also solved the isomorphism problem in \cite{GH} for some multi-parameter quantized Weyl algebras via improving a result of Rigal \cite{Rigal}. These multi-parameter quantized Weyl algebras can be regarded as generalized Weyl algebras of higher ranks. One can always define new families of generalized Weyl algebras starting with different $R$ and $\sigma$. Thus it remains an interesting question to study the automorphism group and solve the isomorphism problem for newly defined families of generalized Weyl algebras. 

In this paper, we will study a new family of generalized Weyl algebras $\A$ and their polynomial extensions. The first task of this paper is to study the isomorphism classification for this family of generalized Weyl algebras $\A$ when none of $p$ and $q$ is a root of unity. Note that the algebra $\mathcal{A}_{p}(1,1, \K_{q}[s,t])$ is isomorphic to the tensor product of the rank-one quantized Weyl algebra $\mathcal{A}_{p}(\K)$ and the quantum plane $\K_{q}[s,t]$. When $p\neq 1$, the automorphism group of $\mathcal{A}_{p}(\K)$ was determined by Alev and Dumas in \cite{AD}. To classify this family of generalized Weyl algebras $\A$ up to algebra isomorphisms, we will need to first determine their automorphism groups. We will first show that the generalized Weyl algebra $\A$ has a simple localization in the case where none of $p$ and $q$ is a root of unity. As a result, we are able to completely determine the height-one prime ideals and the center for $\A$, and prove that $\A$ is cancellative. As an application, we will further determine the automorphism group for $\A$ and solve the isomorphism problem. In addition, we will determine the automorphism group and solve the isomorphism problem for their polynomial extensions. In particular, we will show that each $\K-$algebra automorphism of the polynomial extension $\A[w]$ is triangular.

Dixmier asked the question in \cite{Dixmier} whether each $\K-$algebra endomorphism of the Weyl algebra $A_{n}(\K) $ is actually an automorphism when the base field $\K$ is of zero characteristic. Dixmier's problem has been also referred to as the Dixmier conjecture. It was proved in \cite{BK, Tsuchimoto} that the Dixmier conjecture is stably equivalent to the Jacobian conjecture \cite{Keller}. Recently, there has been some interest in establishing a quantum analogue of the Dixmier conjecture. Speaking of a quantum analogue of the Dixmier conjecture, a `naive' target algebra would be the quantized Weyl algebra $\mathcal{A}_{p}(\K)$. However, it was proved by Backelin in \cite{Backelin} that the $\K-$algebra endomorphisms of $\mathcal{A}_{p}(\K)$ may not necessarily be bijective. We also note that the Weyl algebra $A_{n}(\K)$ is a simple $\K-$algebra when the base field $\K$ is of zero characteristic. Thus a correct version of the quantum Dixmier conjecture would be asking whether each algebra endomorphism of a given simple `quantum' algebra is indeed bijective. It was proved by Richard that each algebra endomorphism of a simple quantum torus is an algebra automorphism in \cite{Richard}. It was also proved that a quantum analogue of the Dixmier conjecture holds for some quantum generalized Weyl algebras in \cite{KL, KL1}. Similar results were established for the simple localizations of some two-parameter down-up algebras in \cite{Tang} and some simple localizations of the multiparameter quantized Weyl algebras in \cite{Tang1}. Since the generalized Weyl algebra $\A$ has a simple localization, it is natural to ask whether a quantum analogue of the Dixmier conjecture holds for the simple localization. The second task of this paper is to establish such a quantum analogue of the Dixmier conjecture for the simple localization $(\mathcal{A}_{p}(1, 1, \K_{q}[s, t]))_{\mathbb{S}}$. 

Alev and Chamarie determined the automorphism groups for the quantum plane $\K_{q}[s,t]$ and their tensor products in \cite{AC} in the case where $q\neq 1$. They also determined the automorphism group for the uni-parameter quantum polynomial algebras $\K_{q}[t_{1}, \cdots, t_{n}]$ when $q$ is not a root of unity. The automorphisms of multi-parameter quantum polynomial algebras $\K_{Q}[t_{1}, \cdots, t_{n}]$ were determined by Osborn and Passman in \cite{OP} in the generic case. Recently, using the methods of discriminants and Mod-p, the automorphism groups have been investigated for the quantum polynomial algebras $\K_{Q}[t_{1}, \cdots, t_{n}]$ and the algebra $\mathcal{A}_{p_{1}}(\K)_{\K}\otimes \cdots \otimes_{\K} \mathcal{A}_{p_{n}}(\K)$ in \cite{CPWZ1, CPWZ2}. Indeed, the method of discriminant applies to many families of PI-algebras, which have dominating discriminants \cite{CPWZ1, CPWZ2}. A method of Mod-p was also introduced in \cite{CPWZ1} to deal with some non-PI algebras. The third task of the paper is to determine the automorphism group for the generalized Weyl algebra $\mathcal{A}_{p}(1, 1, \K_{q}[s,t])$ and its polynomial extension in the case where $p\neq 1$ and $q\neq 1$, using the methods of discriminant and Mod-p as introduced in \cite{CPWZ1, CPWZ2}.

The paper is organized as follows. In Section $1$, we introduce a family of generalized Weyl algebras $\A$ and establish some of their basic properties. In Section $2$, we study the algebra automorphisms and solve the isomorphism problem for the generalized Weyl algebra $\A$ and its polynomial extension $\A[w]$. In Section $3$, we establish a quantum analogue of the Dixmier conjecture for the simple localization $(\mathcal{A}_{p}(1, 1, \K_{q}[s, t]))_{\mathbb{S}}$ and determine its automorphism group. In Section $4$, we determine the automorphism group for the algebra $\mathcal{A}_{p}(1, 1, \K_{q}[s,t])$ and its polynomial extension in the case where $p\neq 1$ and $q\neq 1$.

\section{A family of generalized Weyl algebras}
In this section, we introduce a family of generalized Weyl algebras $\A$ defined over the base ring $R=\K[z]_{q}[s, t]$, which is the $\K-$algebra generated by $z, s, t$ subject to the relation: $zs=sz, zt=tz, st=qts$ where $q\in K^{\ast}$. The algebra $\A$ can be regarded as the extension of the rank-one quantized Weyl algebra $\mathcal{A}_{p}(\K)$ by the quantum plane $\K_{q}[s,t]$. In particular, the algebra $\mathcal{A}_{p}(1,1, \K_{q}[s,t])$ is isomorphic to the tensor product of $\mathcal{A}_{p}(\K)$ and $\K_{q}[s,t]$. We will first introduce these algebras as a family of algebras defined in terms of generators and relations. Then we will realize these algebras as generalized Weyl algebras and establish some basic properties for them.

\begin{defn}
Let $\K$ be a field and $p, q, \lambda, \mu \in \K^{\ast}$. We define the extended quantized Weyl algebra $\A$ to be the $\K-$algebra generated by the generators $x,y, s, t$ subject to the following relations:
\begin{eqnarray*}
xy-pyx=1,\quad\quad st=qts,\\
sx=\lambda xs,\quad\quad sy=\lambda^{-1} ys,\\
tx=\mu xt, \quad\quad ty=\mu^{-1} yt.
\end{eqnarray*}
\end{defn}

Let us set $z=xy-yz$ and it is straightforward to verify the following identities:
\[
xz=pzx,\quad yz=p^{-1} zy,\quad sz=zs, \quad tz=zt.
\]
As a result, we have the following obvious proposition.
\begin{prop}
The algebra $\A$ is an iterated skew polynomial ring in the variables $s, t, x, y$. It is a noetherian domain of Gelfand-Kirillov dimension $4$ and it has a $\K-$basis given as follows:
\[
\{x^{l}y^{m}s^{n} t^{o}\mid l, m, n, o \in \Z_{\geq 0}\}
\]
\end{prop}
{\bf Proof:} It is obvious that we can present $\A$ as an iterated Ore extension in the variables $s, t, x, y$ as follows:
\[
\A=\K[s][t;\sigma_{1}][x;\sigma_{2}][y;\sigma_{3};\delta]
\]
where the automorphisms $\sigma_{1}, \sigma_{2}, \sigma_{3}$ and the derivation $\delta$ are defined as follows:
\[
\sigma_{1}(s)=q^{-1} s,\quad \sigma_{2}(s)=\lambda^{-1} s, \quad \sigma_{2}(t)=\mu^{-1} t
\]
and 
\[
\sigma_{3}(s)=\lambda s,\quad \sigma_{3}(t)=\mu t,\quad \sigma_{3}(x)=p^{-1}
\]
and
\[
\delta(s)=\delta(t)=0,\quad \delta(x)=-p^{-1}.
\]
As a result, we know that $\A$ is a noetherian domain of Gelfand-Kirillov dimension $4$ and has the above $\K-$basis as desired.
\qed 

Now we are going to realize $\A$ as a generalized Weyl algebra. First of all, let us recall the definition of generalized Weyl algebras from \cite{Bavula}. Let $R$ be a ring and $\sigma=(\sigma_{1}, \cdots, \sigma_{n})$ be a set of commuting ring automorphisms of $R$, and $a=(a_{1}, \cdots, a_{n})$ be set of (non-zero) elements in the center of $R$ such that $\sigma_{i} (a_{j})=a_{j}$ for $i\neq j$. The {\it generalized Weyl algebra} $A=R(\sigma, a)$ of degree $n$ with a base ring $R$ is the ring generated by $R$ and $2n$ variables $X_{1}^{\pm}, \cdots, X_{n}^{\pm}$ subject to the relations:
\[
X_{i}^{-} X_{i}^{+}=a_{i}, \quad X_{i}^{+}X_{i}^{-}=\sigma_{i}(a_{i}), \quad X_{i}^{\pm } d =\sigma_{i}^{\pm 1} (d) X_{i}^{\pm }\quad \forall d \in D
\]
and
\[
[X_{i}^{-}, X_{j}^{-}]= [X_{i}^{+}, X_{j}^{+}]=[X_{i}^{+}, X_{i}^{ -}]=0,\quad \forall i\neq j
\]
where $[x, y]:=xy-yx$. If $R$ is a noetherian domain, then $A$ is a noetherian domain. For any vector $k=(k_{1}, \cdots, k_{n})\in \Z^{n}$, let $v_{k}=v_{k1}(1)\cdots v_{kn}(n)$, where for any $1\leq i\leq n$ and $m\geq 0$, $v_{\pm m} (i)=(X_{i}^{\pm })^{m}$ and $v_{0}(i)=1$. It follows from that the definition that 
\[
A=\oplus_{k\in \Z^{n}} A_{k}
\]
is a $\Z^{n}-$graded algebra with $A_{k}=Rv_{k}$. For more details about generalized Weyl algebras, we refer the readers to \cite{Bavula} and the references therein. 

Let $R$ be the subalgebra of $\A$ generated by $z, s, t$. Then $R$ can be regarded a quantum polynomial algebra in two variables over the base ring $\K[z]$. We will sometimes denote $R$ by $\K[z]_{q}[s, t]$. We have the following proposition.
\begin{prop}
The algebra $\A$ is a generalized Weyl algebra over the base algebra $R=\K[z]_{q}[s, t]$. As a result, the algebra $\A$ also has a $\K-$basis given as follows:
\[
\{z^{a}s^{b}t^{c}x^{d}, \, z^{a}s^{b}t^{c}y^{e+1}\mid a, b, c, d, e \in \Z_{\geq 0}\}.
\]
\end{prop}
{\bf Proof:} Note that in $R=\K[z]_{q} [s, t]$, we have the following identities:
\[
xz=pzx,\, yz=p^{-1}yz, \, xs=\lambda^{-1} sx,\, xt=\mu^{-1} tx,\, ys=\lambda sy, \, yt=\mu ty.
\]
Let $\sigma$ be the $\K-$algebra automorphism of $R$ defined as follows:
\[
\sigma(s)=\lambda^{-1} s, \quad \sigma(t)=\mu^{-1} t,\quad \sigma(z)=pz.
\]
Then it is easy to verify that we have the following:
\[
yx=\frac{z-1}{p-1},\quad xy=\frac{pz-1}{p-1}.
\]
As a result, we have that $\A=R(\sigma, \frac{z-1}{p-1})$, which is a generalized Weyl algebra over the base algebra $R=\K[z]_{q}[s, t]$. 
\qed

It is well-known that the rank-one quantized Weyl algebra $\mathcal{A}_{p}(\K)$ has a simple localization with respect to the Ore set generated by $z:=xy-yx$, when $p$ is not a root of unity \cite{Jordan}. The quantum torus $\K_{q}[s^{\pm 1}, t^{\pm 1}]$ is a simple algebra when $q$ is not a root of unity\cite{MP}. Next we show that the algebra $\A$ has a simple localization when none of $p,q$ is not a root of unity. Let us set $\mathbb{S}=\{z^{i}s^{j}t^{k}\mid i, j, k \in \Z_{\geq 0}\}$. Then the set $\mathbb{S}$ is an Ore set for the algebra $\A$. We can localize the algebra $\A$ with respect to the set $\mathbb{S}$ and we denote the localization by $\A_{\mathbb{S}}$.
\begin{thm}
If none of $p, q$ is a root of unity, then the localization $(\A)_{\mathbb{S}}$ is a simple algebra.
\end{thm}
{\bf Proof:} (We follow the proof used by Jordan in \cite{Jordan}). For simplicity, we will denote the algebra $\A$ by $A$ and the localization $(\A)_{\mathbb{S}}$ by $B$. Let us further denote the localization of the algebra $B$ with respect to the Ore set $\mathbb{Y}:=\{y^{i}\mid i\in \mathbb{Z}_{\geq 0}\}$ by $C$. Then it is easy to verify that the algebra $C$ is the quantum torus (or quantum Laurent polynomial algebra) generated by $4$ generators $y, z, s, t$ subject to the following relations:
\[
yz=p^{-1} zy,\quad ys=\lambda sy, \quad yt=\mu ty, \quad zs=sz,\quad zt=tz, \quad st=qts.
\]

Let us place the generators in the order $z, y, s, t$. Then the algebra $C$ is a McConnell-Pettit algebra \cite{MP, MR} generated by the variables $z, y, s, t$ subject to the following commuting relation matrix:
\[
M=\left[
\begin{array}{lccr}
1& p & 1 & 1\\
p^{-1}&1&\lambda&\mu\\
1&\lambda^{-1} &1& q\\
1 & \mu^{-1} & q^{-1} & 1
\end{array}
\right].
\]
As a result, the algebra $C$ is a simple algebra due to the fact that none of $p$ and $q$ is a root of unity \cite{MP}.

Next we need to show that the algebra $B$ is a simple algebra. Let $I$ be a non-zero two-sided ideal of $B$. Then the localization of $I$ with respect to the Ore set $\mathbb{Y}$ is a non-zero two-sided ideal of the algebra $C$. Since $C$ is a simple algebra, we have $I_{\mathbb{Y}}=C$. Thus $I$ must contain an element of the form $y^{j}$ for some $j\geq 0$. Suppose $y^{d}$ is the smallest power of $y$ contained in the ideal $I$ with $d\geq 1$. Note that it is easy to check that we have the following:
\[
xy^{d}-p^{d}y^{d}x=\frac{p^{d}-1}{p-1} y^{d-1}.
\]
Since $y^{d} \in I$ and $I$ is a two-sided ideal of $B$, we have that $\frac{q^{d}-1}{p-1}y^{d-1}\in I$. Sine $p$ is not a root of unity, we have that $y^{d-1} \in I$, which is a contradiction to the assumption that $d$ is the smallest. So we have that $d=0$, which implies that $1\in I$. Thus we have $I=B$. So we have proved that the algebra $B$ is indeed a simple algebra.
\qed

\begin{cor}
The height-one prime ideals of $\A$ are exactly the following ideals:
\[
\{(z),\, (s),\, (t)\}.
\]
\end{cor}
{\bf Proof:} Firt of all, it is obvious that the ideals $(z), (s)$, and $(t)$ are completely prime ideals of $A=\A$ because the quotient algebras $A/(z), A/(s)$ and $A/(t)$ are all domains. Let $P$ be a non-zero height-one prime ideal of $A$. Since $A_{\mathbb{S}}$ is a simple algebra, we have $P_{\mathbb{S}}=A_{\mathbb{S}}$. Thus $P$ must contain an element of the form $z^{i}s^{j}t^{k}$ for some $i, j, k$ not all zero. Since $P$ is a prime ideal and $z, s, t$ are normal elements, $P$ must contain either $z$ or $s$ or $t$. Thus we have that either $(z)\subseteq P$ or $(s) \subseteq P$ or $(t)\subseteq P$. Since $P$ is of height one, we have $P=(z)$ or $P=(s)$ or $P=(t)$. Using the non-commutative Principal Ideal Theorem \cite{MR}, we know that the prime ideals $(z), (s)$ and $(t)$ have height one. Thus we have proved the corollary.
\qed

The study of Zariski cancellation problem for non-commutative algebras has recently been initiated by Bell and Zhang in \cite{BZ}. We now prove that the algebra $\A$ is cancellative. First of all, we need to recall the following definition from \cite{BZ}. 

\begin{defn} Let A be any $\K-$algebra.
\begin{enumerate}
\item We call $A$ cancellative if $A[w]\cong B[w]$ for some $\K-$algebra B implies that $A\cong B$ as a $\K-$algebra.\\

\item We call $A$ strongly cancellative if, for any $d \geq 1$, the isomorphism $A[w_{1}, \cdots, w_{d}] \cong B[w_{1}, \cdots , w_{d}$]
for some $\K-$algebra $B$ implies that $A\cong B$ as a $\K-$algebra.\\

\item We call $A$ universally cancellative if, for any finitely generated commutative
$\K-$algebra and domain $R$ such that $R/I = \K$ for some ideal $I \subset R$ and any $\K-$algebra $B$, $A\otimes R \cong B\otimes R$ implies that $A\cong B$ as a $\K-$algebra.
\end{enumerate}
\end{defn}
\begin{thm}
Assume that none of $p$ and $q$ is a root of unity. Then the center of $\A$ is $\K$. Furthermore, $\A$ is universally cancellative; and thus it is strongly cancellative and cancellative.
\end{thm}
{\bf Proof:} Note that the localization $(\A)_{\mathbb{S}}$ is a simple $\K-$algebra and no element in $\mathbb{S}$ is a central element of $\A$ when none of $p$ and $q$ is a root of unity. Thus the center of $\A$ is $\K$. As a result, $\A$ is a finitely generated domain of Gelfand-Kirillov dimension $4$ and $\A$ is not commutative. Now the result follows from {\bf Theorem 0.5} in \cite{BZ}.
\qed

\begin{rem}
Note that the algebra $\mathcal{A}_{p}(1,1, \K_{q}[s, t])$ is the tensor product of the rank-one quantized Weyl algebra $\mathcal{A}_{p}(\K)$ and the quantum plane $\K_{q}[s,t]$. It was proved in \cite{BZ} that the algebra $\mathcal{A}_{p}( 1, 1, \K_{q}[s, t])$ is cancellative when $p\neq 1$ and $q\neq 1$.
\end{rem}

\section{Automorphisms and Isomorphisms for $\A$}

In this section, we study the automorphism group of the algebra $\A$ and its polynomial extension $\A[w]$ under the condition that none of $p$ and $q$ is a root of unity. As an application, we will solve the isomorphism problem. 

To proceed, we first establish a lemma, whose proof mimics the ones used in \cite{Rigal} and \cite{GH}. And we will not repeat the details of the proof here. Let us set $z_{0}=1, z_{1}=z, z_{2}=s, z_{3}=t$. We have the following lemma.
\begin{lem}
Let $a,b \in \A \backslash \K$ such that 
\[
ab=\sum_{j=0}^{3} \alpha_{j} z_{j}
\]
with $\alpha_{j}\in \K$ and $\alpha_{1}\neq 0$. Then there exist $\lambda, \mu \in \K^{\ast}$ such that 
\[
a=\lambda y,\quad \quad b=\mu x
\]
or 
\[
a=\lambda x,\quad\quad b=\mu y.
\]
\end{lem}
\qed

\begin{thm}
Assume that none of $p,q$ is a root of unity. Let $\varphi$ be a $\K-$algebra automorphism of $\A$. Then we have the following:
\[
\varphi(x)=\alpha x, \quad \varphi(y)=\alpha^{-1} y, \quad \varphi(s)=\beta s,\quad \varphi(t)=\gamma t.
\]
for some $\alpha, \beta, \gamma\in \K^{\ast}$.
\end{thm}

{\bf Proof:} Since $\varphi$ is a $\K-$algebra automorphism of $\A$, it permutes all the height-one prime ideals $(z), (s)$, and $(t)$ of the algebra $\A$. As a result, $\varphi$ has to permute the elements $z, s$, and $t$ up to nonzero scalars. Thus we can restrict the automorphism $\varphi$ to the subalgebra $R$ of $\A$ generated by $z, s, t$. It is easy to see that the restriction of $\varphi$ to the subalgebra $R$ is also a $\K-$algebra automorphism of $R$. Since the center of $R$ is $\K[z]$, we shall have that $\sigma(z)=a z$ for some $a \in \K^{\ast}$. As a result, we know that $\varphi$ permutes $s, t$ up to non-zero scalars. Since $st=qts$ and $q$ is not a root of unity, we shall further have that $\sigma(s)=\beta s$ and $\sigma(t)=\gamma t$ for some $\beta, \gamma \in \K^{\ast}$. 

Since $A$ is a generalized Weyl algebra over $R$, it has a $\K-$ basis given as follows:
\[
\{z^{l}s^{m}t^{n}x^{i}, \, z^{l}s^{m}t^{n}y^{j+1}\mid l, m,n, i, j \in \mathbb{Z}_{\geq 0}\}.
\] 
Therefore, we can express the elements $\varphi(x)$ and $\varphi(y)$ in terms of this new basis as follows:
\[
\varphi(x)=\sum_{i\geq 0} f_{i}(z, s, t) x^{i}+\sum_{j\geq 1} g_{j}(z, s, t) y^{j}
\]
where $f_{i}(z, s, t), g_{j}(z, s, t) \in R$. Since $xz=pzx$, we have $\varphi(x)\varphi(z)=p\varphi(z)\varphi(x)$. As a result, we have the following:
\[
\varphi(x)az=p az\varphi(x)
\]
which implies that $p^{i}=p$ for such $i$ with $f_{i}\neq 0$ and $p^{-j}=p$ for such $j$ with $g_{j}(z, s, t)\neq 0$. Since $p$ is not a root of unity, we have $i=1$ and $j=-1$. Since $j\geq 0$, we have $g_{j}(z, s, t)=0$ and $f_{i}(z, s, t)=0$ for $i\neq 1$. So we have the following:
\[
\varphi(x)=f(z, s, t) x
\]
for some $f(z, s, t) \in R$. Similarly, we shall have that $\varphi(y)=g(z, s, t) y$ for some $g(z, s, t) \in R$. Since $yx=\frac{z-1}{p-1}$ and $\varphi(z)=az$, we have the following:
\begin{eqnarray*}
\frac{az-1}{p-1} &=&\varphi(y)\varphi(x)\\
&=&g(z, s, t) yf(z, s, t) x\\
&=&f(p^{-1}z, \lambda s, \mu t)g(z, s, t)yx\\
&=& \frac{f(p^{-1}z, \lambda s, \mu t)g(z, s, t)(z-1)}{p-1}.
\end{eqnarray*}
Thus, we have that $f(p^{-1}z, \lambda s, \mu t)g(z, s, t)=1$ and $a=1$. As a result, both $f(z, s, t)$ and $g(z, s, t)$ are in $\K^{\ast}$. Therefore, we have the following:
\[
\varphi(x)=\alpha x,\quad \varphi(y)=\alpha^{-1}y, \quad \varphi (z)=z
\]
for some $\alpha \in \K^{\ast}$. So we have completed the proof.
\qed

For any scalars $\alpha, \beta, \gamma \in \K^{\ast}$, it is easy to check the following mapping
\[
\varphi(x)=\alpha x, \quad \varphi(y)=\alpha^{-1} y, \quad \varphi(t_{1})=\beta t_{1},\quad \varphi(t_{2})=\gamma t_{2}
\]
defines a $\K-$algebra automorphism of $\A$. And we denote this automorphism by $\varphi_{\alpha, \beta, \gamma}$. Let ${\rm Aut}_{\K}(A)$ denote the group of all $\K-$algebra automorphisms of $\A$. Then we have the following result.
\begin{thm} The automorphism group of $\A$ is given as follows:
\[
{\rm Aut}_{\K}(\A)=\{\varphi_{\alpha, \beta, \gamma}\mid \alpha, \beta, \gamma \in \K^{\ast}\}\cong (\K^{\ast})^{3}.
\]
\end{thm}
{\bf Proof:} The result follows from the fact that any $\K-$algebra automorphism of $\A$ is of the form $\varphi_{\alpha, \beta,\gamma}$, which is completely determined by the scalars $\alpha, \beta, \gamma$, and $\varphi_{\alpha, \beta, \gamma}\circ \varphi_{\alpha^{\prime}, \beta^{\prime}, \gamma^{\prime}}=\varphi_{\alpha\alpha^{\prime}, \beta\beta^{\prime}, \gamma\gamma^{\prime}}$, and $\varphi_{1, 1, 1}=id$, and $\varphi_{\alpha, \beta, \gamma}^{-1}=\varphi_{\alpha^{-1}, \beta^{-1}, \gamma^{-1}}$. 
\qed

\begin{thm}
Assmue that none of $p_{1}, p_{2}$ and $q_{1}, q_{2}$ is a root of unity. Let $A_{1}=\D$ and $A_{2}=\E$ be two extended quantized Weyl algebras generated by $x_{1}, y_{1}, s_{1}, t_{1}$ and $x_{2}, y_{2}, s_{2}, t_{2}$ respectively. Then $A_{1}$ is isomorphic to $A_{2}$ as a $\K-$algebra if and only if 
\[
p_{1}=p_{2}, \quad q_{1}=q_{2},\quad \lambda_{1}=\lambda_{2},\quad \mu_{1}=\mu_{2}
\]
or
\[
p_{1}=p_{2}^{-1}, \quad q_{1}=q_{2}, \quad \lambda_{1}=\lambda_{2}^{-1}, \quad \mu_{1}=\mu_{2}^{-1}
\]
or
\[
p_{1}=p_{2},\quad q_{1}=q_{2}^{-1}, \quad \lambda_{1}=\mu_{2},\quad \lambda_{2}=\mu_{1}
\]
or
\[
p_{1}=p_{2}^{-1}, \quad q_{1}=q_{2}^{-1}, \quad \lambda_{1}=\mu_{2}^{-1}, \quad \lambda_{2}=\mu_{1}^{-1}.
\]
\end{thm}
{\bf Proof:} Let $\psi$ be a $\K-$algebra isomorphism from $A_{1}$ to $A_{2}$. Then $\psi$ sends all the height one prime ideals $(z_{1}), (s_{1}), (t_{1})$ of $A_{1}$ to the height one prime ideals $(z_{2}), (s_{2}), (t_{2})$ of $A_{2}$. As a result, $\psi$ sends the elements $z_{1}, s_{1}, t_{1}$ to scalar multiples of $z_{2}, s_{2}, t_{2}$. Thus we can restrict the isomorphism $\psi$ to an isomorphism from the subalgebra $R_{1}=\K[z_{1}]_{q_{1}}[s_{1}, t_{1}]$ to the subalgebra $R_{2}=\K[z_{2}]_{q_{2}}[s_{2}, t_{2}]$. Since $z_{1}$ is central in $R_{1}$ and $z_{2}$ is central in $R_{2}$, we shall have that $\psi(z_{1})=a z_{2}$ for some $a\in \K^{\ast}$. Moreover, we shall have that $q_{1}=q_{2},\psi(s_{1})=\beta s_{2},\psi(t_{1})=\gamma t_{2} $ or $\psi(s_{1})=\beta t_{2}, \psi(t_{1})=\gamma s_{2}, q_{1}=q_{2}^{-1}$. Using an analogue of {\bf Lemma 2.1}, we shall have that either $\sigma(x_{1})=\alpha x_{2}, \psi(y_{1})=\alpha^{-1} y_{2}, p_{1}=p_{2}$ or $\psi(x_{1})=\alpha y_{2}, \psi(y_{1})=-\alpha^{-1} p_{1}^{-1}x_{2}, p_{1}=p_{2}^{-1}$. The conditions on the parameters $\lambda_{1}, \lambda_{2}, \mu_{1}, \mu_{2}$ also follow. Conversely, if one of the conditions is satisified, it is straightforward to check that we can define a $\K-$algebra isomorphism from $A_{1}$ to $A_{2}$. So we have completed the proof.
\qed

We now study the $\K-$algebra automorphisms of the polynomial extension of $\A$. First of all, we need to recall a notion. Let $A$ be any $\K-$algebra. A $\K-$algebra automorphism $h$ of the polynomial extension $A[w]=A\otimes_{\K} K[w]$ is said to be {\it triangular} if there is a $g\in {\rm Aut}_{\K} (A)$ and $c\in \K^{\ast}$ and $r$ in the center of $A$ such that 
\[
h(w)=cw+r
\]
and
\[
h(u)=g(u)\in A
\]
for any $u\in A$. We have following theorem about the $\K-$algebra automorphisms of $\A[w]$ and locally nilpotent $\K-$derivations of $\A$.
\begin{thm}
Suppose that none of $p$ and $q$ is a root of unity. Any $\K-$algebra automorphism $\varphi$ of $\A[w]$ is given as follows
\[
\varphi(x)=\alpha x, \quad \varphi(y)=\alpha^{-1} y, \quad \varphi(s)=\beta s,\quad \varphi(t)=\gamma t
\]
for some $\alpha, \beta, \gamma\in \K^{\ast}$ and
\[
\varphi(w)=aw+b
\]
for some $a\in \K^{\ast}$ and $b\in \K$. As a result, any $\K-$algebra automorphism $\varphi$ of $\A$ is triangular. If $\Z\subset \K$, then any locally nilpotent $\K-$derivation of $\A$ is zero.
\end{thm}
{\bf Proof:} Let $\mathbb{W}=\K[w]-\{0\}$. Then $\mathbb{W}$ is an Ore set for the algebra $\A[w]$. It is obvious that the corresponding localization $(\A[w])_{\mathbb{W}}$ is isomorphic to the generalized Weyl algebra $\mathcal{A}_{p}(\lambda, \mu, \K(w)_{q}[s,t])$ over the base ring $\K(w)_{q}[s,t]$. As a result, the center of $\A[w]$ is $\K[w]$. Let $\varphi$ be a $\K-$algebra automorphism of $\A[w]$. Then $\varphi$ preserves the subalgebra $\K[w]$ and thus it can be restricted to a $\K-$algebra automorphism of $\K[w]$. As a result, we have the following:
\[
\varphi(w)=aw+b
\]
for some $a\in \K^{\ast}$ and $b\in \K$. Let $\varphi_{0}$ be the $\K-$algebra automorphism of $\A[w]$ defined as follows:
\[
\varphi_{0}(x)=x, \,\varphi_{0}(y)=y,\, \varphi_{0}(s)=s, \,\varphi_{0}(t)=t,\, \varphi_{0}(w)=a^{-1}w-a^{-1}b.
\]
Now the composition $\varphi_{0}\circ \varphi$ of $\varphi_{0}$ and $\varphi$ is a $\K-$algebra automorphism of $\A[w]$ which fixes the variable $w$. Thus we can regard the composition $\psi\circ \varphi$ as a $\K[w]-$algebra automorphism of $\A[w]$, which can be further regarded as $\K(w)-$algebra automorphism of $(\A)_{\mathbb{W}}=\mathcal{A}_{p}(\lambda, \mu, \K(w)_{q}[s,t])$. By {\bf Theorem 2.1}, we have the following:
\[
\varphi_{0}\circ \varphi(x)=\alpha x, \quad \varphi_{0}\circ \varphi(y)=\alpha^{-1} y, \quad \varphi_{0}\circ \varphi(s)=\beta s,\quad \varphi_{0}\circ \varphi(t)=\gamma t.
\]
for some $\alpha, \beta, \gamma\in \K(w)^{\ast}$. As a result, we have the following
\[
\varphi(x)=\alpha x, \quad \varphi(y)=\alpha^{-1} y, \quad \varphi(s)=\beta s,\quad \varphi(t)=\gamma t, \varphi(w)=aw+b.
\]
for some $\alpha, \beta, \gamma\in \K(w)^{\ast}$ and $a\in \K^{\ast}$ and $b\in \K$. Since $\varphi$ is a $\K-$algebra automorphism of $\A[w]$, we shall have $\alpha, \beta, \gamma \in \K[w]$. Since $\varphi$ has an inverse, we shall further have that $\alpha, \beta, \gamma \in \K^{\ast}$ as desired. Therefore, any $\K-$algebra automorphism $\varphi$ of $\A[w]$ is triangular. If $\Z\subset \K$, we can further prove that any locally nilpotent $\K-$derivation of $\A$ is zero using the argument as outlined in \cite{CPWZ1}.
\qed

Moreover, we have the following result on the isomorphism problem for the polynomial extensions of the generalized Weyl algebras $\A$. 
\begin{thm}
Assume that none of $p_{1}, p_{2}$ and $q_{1}, q_{2}$ is a root of unity. Then $\D[w]$ is isomorphic to $\E[w]$ as a $\K-$algebra if and only if 
\[
p_{1}=p_{2}, \quad q_{1}=q_{2},\quad \lambda_{1}=\lambda_{2},\quad \mu_{1}=\mu_{2}
\]
or
\[
p_{1}=p_{2}^{-1}, \quad q_{1}=q_{2}, \lambda_{1}=\lambda_{2}^{-1}, \quad \mu_{1}=\mu_{2}^{-1}
\]
or
\[
p_{1}=p_{2}, \quad q_{1}=q_{2}^{-1}, \quad \lambda_{1}=\mu_{2}, \quad \lambda_{2}=\mu_{1}
\]
or
\[
p_{1}=p_{2}^{-1}, \quad q_{1}=q_{2}^{-1}, \quad \lambda_{1}=\mu_{2}^{-1},\quad \lambda_{2}=\mu_{1}^{-1}.
\]
\end{thm}
{\bf Proof:} Since the algebra $\D$ is cancellative, we know that $\D[w]\cong \E[w]$ if and only if $\D\cong \E$. Now the result follows directly from {\bf Theorem 2.3}. 
\qed

\section{Endomorphisms for the simple localization $(\mathcal{A}_{p}(1, 1, \K_{q}[s,t]))_{\mathbb{S}}$}
When none of the parameters $p$ and $q$ is a root of unity, the localization $(\mathcal{A}_{p}(1, 1, \K_{q}[s,t]))_{\mathbb{S}}$ is simple. Therefore, one may ask the question whether a quantum analogue of the Dixmier conjecture holds for $(\mathcal{A}_{p, 1, 1}(\K_{q}[s,t]))_{\mathbb{S}}$. In this section, we will show that each $\K-$algebra endomorphism of $(\mathcal{A}_{p, 1, 1}(\K_{q}[s,t]))_{\mathbb{S}}$ is indeed an algebra automorphism. Moreover, we will completely describe the automorphism group for this simple localization. 

In the rest of this section, we will denote the algebra $\mathcal{A}_{p}(1,1, \K_{q}[s, t])$ by $A$ and denote its localization with respect to $\mathbb{S}$ by $A_{\mathbb{S}}$. First of all, it is obvious that the invertible elements of the algebra $A_{\mathbb{S}}$ are exactly of the form $az^{l}s^{m}t^{n}$ where $a\in \K^{\ast}$ and $l, m, n \in \Z$. As a result, we have the following lemma.

\begin{lem}
Let $\varphi$ be a $\K-$algebra endomorphism of $A_{\mathbb{S}}$, then we have the following:
\[
\varphi(x)=g(z, s, t) x, \quad \varphi(y)=h(z, s, t) y,\quad \varphi(z)=\alpha z\]
or 
\[
\varphi(x)=g(z, s, t) y, \quad \varphi(y)=h(z, s, t) x, \quad \varphi(z)=\alpha z^{-1} 
\]
and 
\[
\varphi(s)=\beta z^{m} s^{c} t^{d}, \quad \varphi(t)=\gamma z^{n} s^{e} t^{f}
\]
for some $g(z, s, t), h(z, s, t) \in \K[z^{\pm 1}]_{q}[s^{\pm 1}, t^{\pm 1}]$ and $\alpha, \beta, \gamma \in \K^{\ast}$ such that $cf-de=1$.
\end{lem}

{\bf Proof:} Since $\varphi$ is a $\K-$algebra endomorphism of $A_{\mathbb{S}}$ and the elements $z,s$ and $t$ are invertible elements of $A_{\mathbb{S}}$, the images $\varphi(z), \varphi(s)$ and $\varphi(t)$ are also invertible elements of $A_{\mathbb{S}}$. As a result, we have the following:
\[
\varphi(z)=\alpha z^{l} s^{a}t^{b}, \quad \varphi(s)=\beta z^{m} s^{c} t^{d}, \quad \varphi(t)=\gamma z^{n} s^{e} t^{f}
\]
where $\alpha, \beta, \gamma \in \K^{\ast}$ and $l, m, n, a, b, c, d \in \Z$. 

Applying $\varphi$ to $zs=sz, zt=tz, st=qts$, we have the following
\[
\varphi(z)\varphi(s)=\varphi(s)\varphi(z),\quad \varphi(z)\varphi(t)=\varphi(t)\varphi(z),\quad \varphi(s)\varphi(t)=q\varphi(t)\varphi(s).
\]
As a result, we have that $q^{ad-bc}=1$, $q^{af-be}=1$, and $q^{cf-de}=q$. Since $q$ is not a root of unity, we have that $ad-bc=0$, $af-be=0$, and $cf-de=1$. Therefore, we know that $a=b=0$ and 
\[
\varphi(z)=\alpha z^{l}, \quad \varphi (s)=\beta z^{m} s^{c} t^{d},\quad \varphi(t)=\gamma z^{n} s^{e} t^{f}.
\]

Let $\varphi(x)=\sum_{i} g_{i}(z, s, t) X^{i}$, where $g_{i}(z, s, t) \in \K[z^{\pm 1}]_{q}[s^{\pm 1}, t^{\pm 1}]$ and $X^{i}=x^{i}$ if $i\geq 0$ and $X^{i}=y^{i}$ if $i<0$. Applying $\varphi$ to $xz=pzx$, we have the following
\[
(\sum_{i} g_{i}(z,s,t) X^{i})(\alpha z^{l})=p (\alpha z^{l})( \sum_{i} g_{i}(z, s,t) X^{i}).
\]
If $g_{i}(z, s, t) \neq 0$, then $li=1$. Thus, we have that $l=\pm 1$ and $i=\mp 1$. Since $z=xy-yx$ and $\varphi(z)\neq 0$, we cannot have either $\varphi(x)=0$ or $\varphi(y)=0$. Indeed, $\varphi(x)=0$ implies that $\varphi(z)=0$, which is a contradiction. As a result, we have that $\varphi(x)=g(z, s, t) x$ and $\varphi(z)=\alpha z$, or $\varphi(x)=g(z, s, t) y$ and $\varphi(z)=\alpha z^{-1}$. Therefore, we shall have that $\varphi(y)=h(z, s, t) y$ if $\varphi(z)=\alpha z$, or $\varphi(y)=h(z, s, t) x$ if $\varphi(z)=\alpha z^{-1}$.
\qed

From now on, we will say that a $\K-$algebra endomorphism $\varphi$ of $A_{\mathbb{S}}$ is of positive type if $\varphi(z)=\alpha z$, and it is of negative type if $\varphi(z)=\alpha z^{-1}$. If $\varphi$ is a $\K-$algebra endomorphism of negative type, then $\varphi^{2}=\varphi\circ \varphi$ is of positive type. We will show each $\K-$algebra endomorphism $\varphi$ of $A_{\mathbb{S}}$ is indeed an automorphism if it is of positive type. As a result, we can conclude that each $\K-$algebra endomorphism $\varphi$ of $A_{\mathbb{S}}$ is an automorphism.

\begin{thm}
Each $\K-$algebra endomorphism $\varphi$ of positive type for the algebra $A_{\mathbb{S}}$ is an automorphism.
\end{thm}
{\bf Proof:} Let $\varphi$ be a $\K-$algebra endomorphism of $A_{\mathbb{S}}$, which is of positive type. By {\bf Lemma 3.1}, we have the following:
\[
\varphi(x)=g(z, ,s, t) x, \quad \varphi(y)=h(z, s, t) y,\quad \varphi(z)=\alpha z\]
and
\[
\varphi(s)=\beta z^{m} s^{c} t^{d}, \quad \varphi(t)=\gamma z^{n} s^{e} t^{f}
\]
for some $g(z, s, t), h(z, s, t) \in \K[z^{\pm 1}]_{q}[s^{\pm 1}, t^{\pm 1}]$ and $\alpha, \beta, \gamma \in \K^{\ast}$ such that $cf-de=1$. Applying $\varphi$ to $xy=\frac{pz-1}{p-1}$, we have the following: 
\begin{eqnarray*}
\varphi(x)\varphi(y)&=&g(z, s, t) x h(z, s, t) y\\
&=&g(z, s, t) h^{\prime}(z, s, t) xy\\
&=&g(z, s, t) h^{\prime}(z, s, t) \frac{pz-1}{p-1}\\
&=&\frac{p \alpha z-1}{p-1}
\end{eqnarray*}
where $h^{\prime}(z, s, t)=h(pz, s, t)$. Thus, we have that $g(z, s, t) h^{\prime}(z, s, t)=1$ and $\alpha=1$. As a result, we have that both $g(z, s, t)$ and $h^{\prime}(z, s, t)$ are invertible. Therefore, both $g(z, s, t)$ and $h(z, s, t)$ are invertible. In particular, we have that
\[
\varphi (x)=\alpha z^{i} s^{j} t^{k} x, \quad \varphi (y)=\alpha^{-1}p^{i}q^{-jk} z^{-i} s^{-j} t^{-k} y, \quad \varphi(z)=z
\]
and 
\[
\varphi(s)=\beta z^{m} s^{c} t^{d}, \quad \varphi(t)=\gamma z^{n} s^{e}t^{f}
\]
such that $p^{m} q^{dj-ck}=1, p^{n} q^{fj-ek}=1$, and $cf-de=1$. Conversely, given any $\alpha, \beta, \gamma \in \K^{\ast}$ and $m, n, c, d, e, f, i, j, k \in \Z$ satisfying the condition: $p^{m} q^{dj-ck}=1, p^{n} q^{fj-ek}=1$ and $cf-de=1$ , we can define a $\K-$algebra endomorphism of $A_{\mathbb{S}}$ and we denote it by $\varphi_{\alpha, \beta, \gamma,i, j, k}^{c, d, e, f, m, n}$.

Let $j^{\prime}=ek-fj$ and $k^{\prime}=dj-ck$ and $m^{\prime}=dn-fm$ and $n^{\prime}=em-cn$ and $c^{\prime}=f$ and $f^{\prime}=c$ and $e^{\prime}=-e$ and $d^{\prime}=-d$. Then it is straightfoward to check that $p^{m^{\prime}} q^{j^{\prime}d^{\prime}-c^{\prime}k^{\prime}}=1, p^{n^{\prime}} q^{j^{\prime}f^{\prime}-k^{\prime}e^{\prime}}=1$, and $c^{\prime}f^{\prime}-d^{\prime}e^{\prime}=1$. Let $\psi=\varphi_{1, 1, 1, 0, j^{\prime}, k^{\prime}}^{c^{\prime}, d^{\prime}, e^{\prime}, f^{\prime}, m^{\prime}, n^{\prime}}$ be the $\K-$algebra endomorphism of $A_{\mathbb{S}}$ defined by the data $\alpha=\beta=\gamma=1, i^{\prime}=0, j^{\prime}, k^{\prime}, c^{\prime}, d^{\prime}, e^{\prime}, f^{\prime}, m^{\prime}$ and $n^{\prime}$. Therefore, we have that $\psi(x)=s^{j^{\prime}} t^{k^{\prime}}x=s^{ek-fj}t^{dj-ck}x$ and
\[
\psi(y)=q^{-j^{\prime}k^{\prime}}s^{-j^{\prime}} t^{-k^{\prime}} y=q^{(dj-ck)(fj-ek)}s^{fj-ek}t^{ck-dj}y,\quad \psi (z)=z
\]
and
\[
\psi(s)=z^{m^{\prime}} s^{c^{\prime}} t^{d^{\prime}}=z^{dn-fm} s^{f} t^{-d},\quad
\psi(t)=z^{n^{\prime}} s^{e^{\prime}} t^{f^{\prime}}=z^{em-cn} s^{-e} t^{c}.
\]
Now the composition $\varphi \circ \psi$ is also a $\K-$algebra endomorphism of $A_{\mathbb{S}}$, which is indeed of positive type. In particuar, it is easy to check that we have the following:
\[
\varphi \circ \psi(x)=\alpha^{\prime} z^{i^{\prime\prime}} x,\quad \varphi\circ \psi(y)=(\alpha^{\prime})^{-1}p^{i^{\prime\prime}}z^{-i^{\prime\prime}}y
\]
and
\[
\varphi \circ \psi(s)=\beta^{\prime} s, \quad \varphi \circ \psi(t)=\gamma^{\prime} t,\quad \varphi\circ \psi (z)=z
\]
where $\alpha^{\prime}, \beta^{\prime},\gamma^{\prime}\in \K^{\ast}$ and $i^{\prime\prime}\in \Z$. It is obvious that $\varphi\circ \psi$ is a $\K-$algebra automorphism of $A_{\mathbb{S}}$. Therefore, $\varphi$ is surjective. Since $\A_{\mathbb{S}}$ is a simple $\K-$algebra, $\varphi$ is also injective. Therefore, $\varphi$ is a $\K-$algebra automorphism of $A_{\mathbb{S}}$ as desired.
\qed

\begin{prop}
Let $\varphi$ be a negative $\K-$algebra automorphism of $A_{\mathbb{S}}$. Then we have the following
\[
\varphi(x)=\alpha p^{-i}z^{-i-1} s^{j}t^{k}y, \quad \varphi(y)=\alpha^{-1} p^{2i} q^{-jk}z^{i} s^{-j} t^{-k} x
\]
and
\[
\varphi(s)=\beta z^{m}s^{c}t^{d}, \quad \varphi(t)=\gamma z^{n} s^{e} t^{f}
\]
for some $i, j, k, m, n, c, d, e, f \in \Z$ and $\alpha, \beta, \gamma \in \K^{\ast}$ such that 
\[
cd-ef=1,\quad p^{m}q^{jd-ck}=1, \quad p^{n}q^{jf-ke}=1.
\] 
Conversely, the above mapping defines a negative $\K-$algebra automorphism of $A_{\mathbb{S}}$.
\end{prop}
{\bf Proof:} First of all, let us define a $\K-$linear mapping $\psi_{-}$ of $A_{\mathbb{S}}$ as follows:
\[
\psi_{-}(x)=y, \, \psi_{-} (y)=-p^{-1} z^{-1} x,\, \psi_{-} (z)=p^{-1} z^{-1},\, \psi_{-}(s)=s, \, \psi_{-} (t)=t.
\]
Then it is easy to see that $\psi_{-}$ is a negative $\K-$algebra automorphism of $A_{\mathbb{S}}$. Indeed, the $\K-$algebra inverse of $\psi$ is defined as follows:
\[
\psi_{-}^{-1}(x)=-z^{-1}y, \, \psi_{-}^{-1}(y)=x,\, \psi_{-}^{-1}(z)=p^{-1}z^{-1},\, \psi(s)=s, \, \psi_{-}(t)=t.
\]
Obviously, $\psi_{-}^{-1}$ is also a negative $\K-$algebra automorphism of $A_{\mathbb{S}}$. Let $\varphi_{-}$ be any negative $\K-$algebra automorphism of $A_{\mathbb{S}}$, then the composition $\psi_{-} \circ \varphi_{-}$ is a positive $\K-$algebra automorphism of $A_{\mathbb{S}}$. As a result, any negative $\K-$algebra automorphism $\varphi_{-}$ is the composition of $\psi_{-}^{-1}$ with a positive $\K-$algebra automorphism $\varphi_{+}$ of $A_{\mathbb{S}}$. Therefore, we have the above description for the negative $\K-$algebra automorphisms of $A_{\mathbb{S}}$.
\qed

Now we give a complete description for the $\K-$algebra automorphism group of $A_{\mathbb{S}}$. Recall that $\psi_{-}$ is the $\K-$algebra automorphism of $A_{\mathbb{S}}$ defined in the proof of {\bf Proposition 3.1}. We denote by $G_{+}$ the subgroup of ${\rm Aut}_{\K}(A_{\mathbb{S}})$, which consists of all the positive $\K-$algebra automorphisms of $A_{\mathbb{S}}$; and let us set $G_{-}=\{\psi_{-}^{-1} \circ \varphi_{+} \mid \varphi_{+} \in G_{+}\}$. Then we have the following result.
\begin{thm}
Any $\K-$algebra automorphism of $A_{\mathbb{S}}$ is either of positive type or the composition of a $\K-$algebra automorphism $\varphi_{+}$ for $A_{\mathbb{S}}$ of postive type with the above defined negative $\K-$algebra automorphism $\psi_{-1}$. In particular, we have the following:
\[
{\rm Aut}_{\K}(A_{\mathbb{S}})=G_{+} \ensuremath{\mathaccent\cdot\cup} G_{-}
\]
and a short exact sequence
\[
1\longrightarrow G_{+}\longrightarrow {\rm Aut}_{\K} (A_{\mathbb{S}})\longrightarrow \Z_{2}\longrightarrow 1.
\]
\end{thm}
{\bf Proof:} Note that the composition of $\psi_{-}^{-1}$ with any positive $\K-$algebra automorphism $\varphi_{+}$ of $A_{\mathbb{S}}$ is a negative $\K-$algebra automorphism for $A_{\mathbb{S}}$. Conversely, the composition of any negative $\K-$algebra automorphism $\varphi_{-}$ of $A_{\mathbb{S}}$ with $\psi$ is a positive $\K-$algebra automorphism for $A_{\mathbb{S}}$. As a result, any negative $\K-$algebra automorphism $\varphi_{-}$ is the composition of $\psi_{-}^{-1}$ with a positive $\K-$algebra automorphism $\varphi_{+}$ of $A_{\mathbb{S}}$. Note that each $\K-$algebra automorphism $\varphi$ of $A_{\mathbb{S}}$ is either of positive type or negative type, but not both. Now the result follows as desired.
\qed

\section{The Automorphism Group of the algebra $\mathcal{A}_{p}(1,1, \K_{q}[s, t])$}
In this section, we will apply the methods of discriminant and Mod-p to determine the $\K-$algebra automorphism group for the algebra $\mathcal{A}_{p}(1, 1, \K_{q}[s,t])$ and its polynomial extension $\mathcal{A}_{p}(1, 1, \K_{q}[s, t])[w]$ in the case where $p\neq 1$ and $q\neq 1$. Note that the algebra $\mathcal{A}_{p}(1, 1, \K_{q}[s,t])$ belongs to the category {\bf Af} as defined in \cite{CPWZ1} when both $p$ and $q$ are roots of unity. For more details about the discriminant and the corresponding terminologies, we refer the readers to \cite{CPWZ1, CPWZ2} and the references there. 

Let $A$ be any filtered $\K-$algebra with a filtration $\{F_{i}A\}_{i\geq 0}$ such that the associated graded algebra ${\rm gr}A$ is generated in degree $1$. A $\K-$algebra automorphism $g$ of $A$ is called {\it affine} if $g(F_{1}A)\subseteq F_{1}A$. First of all, we have the following result.
\begin{thm}
Assume that $p\neq 1$ and $q\neq 1$. Then each $\K-$algebra automorphism of $\mathcal{A}_{p}(1, 1, \K_{q}[s, t])$ is affine. 
\end{thm}
{\bf Proof:} We have a similar proof to the one used for the tensor product of rank-one quantized Weyl algebras in \cite{CPWZ2}. For the reader's convenience, we will produce the details here. Let $Y$ be the subspace $(\K x \oplus \K y) \oplus (\K s\oplus \K t)$. We can regard $\mathcal{A}_{p}(1,1, \K_{q}[s, t])$ as a filtered $\K-$algebra with the standard filtration $F_{n}:=(Y\oplus \K)^{n}$ with degrees of $x, y, s, t$ all equal to $1$. It is easy to see that the associated graded algebra of $\mathcal{A}_{p}(1,1, \K_{q}[s, t])$ is a skew polynomial ring in the variables $x, y, s, t$. Thus $\mathcal{A}_{p}(1,1, \K_{q}[s, t])$ has a monomial basis given as follows
\[
\{ x^{l}y^{m} s^{n} t^{o} \mid l, m, n, o\in \Z_{\geq 0}\}.
\] 
Let $\varphi$ be a $\K-$algebra automorphism of $\mathcal{A}_{p}(1,1, \K_{q}[s, t])$ and write $\varphi (x), \varphi (y), \varphi (s), \varphi(t)$ and $\varphi^{-1}(x), \varphi^{-1}(y), \varphi^{-1}(s), \varphi^{-1} (t)$ as linear combinations of monomials in $x, y, s,t$. Let $\mathcal{K}$ be the $\Z-$subalgebra of the field $\K$ generated by the coefficients of these linear combinations along with their inverses and $p^{\pm 1}$ and $q^{\pm 1}$ and $(p-1)^{\pm 1}$. Let $B$ be the corresponding $\mathcal{K}-$subalgebra of $\mathcal{A}_{p}(1,1, \K_{q}[s, t])$ generated by $x,y, s, t$. Then it is obvious that both $\varphi$ and $\varphi^{-1}$ can be restricted to $B$ as $\mathcal{K}-$algebra automorphisms. $B$ is again an iterated Ore extension of the base field $\mathcal{K}$ and $B$ is free over $\mathcal{K}$. Let $F=\mathcal{K}/\mathfrak{m}$ be any finite quotient field of $\mathcal{K}$. Then $B\otimes_{\mathcal{K}} F$ is also iterated Ore extension, which is isomorphic to the tensor product of the rank-one quantized Weyl algebra and the quantum plane over the finite field $F$. Thus it is in the category {\bf Af} as defined in \cite{CPWZ1, CPWZ2}. As a result, the restriction of $\varphi$ to the subalgebra $B$ is affine. Therefore, $\varphi$ is affine as well, using {\bf Lemma 4.6} in \cite{CPWZ2}.
\qed

For simplicity, we will further denote $\mathcal{A}_{p}(1, 1, \K_{q}[s, t])$ by $B$ and let ${\rm LNDer}(B)$ denote the set of all locally nilpotent $\K-$linear derivations of $B$. Let $B[w]=B\otimes_{\K} \K[w]$ denote the polynomial extension of $B$ by $\K[w]$. We have the following theorem.
\begin{thm}
Assume that $p\neq 1$ and $q\neq 1$. The we have the following result.
\begin{enumerate}

\item Each $\K-$algebra automorphism of $B$ is affine. Moreover, we have the following:
\begin{enumerate}

\item If $p\neq -1$ and $q\neq -1$, then ${\rm Aut}_{\K}(B)=(\K^{\ast})^{3}$. \\

\item If $p=-1$ and $q\neq -1$, then ${\rm Aut}_{\K}(B)=\Z_{2}\ltimes (\K^{\ast})^{3}$. \\

\item If $p\neq -1$ and $q=-1$, then ${\rm Aut}_{\K}(B)=\Z_{2}\ltimes (\K^{\ast})^{3}$. \\

\item If $p=-1$ and $q=-1$, then ${\rm Aut}_{\K}(B)=\Z_{2}\ltimes (\Z_{2}\ltimes (\K^{\ast})^{3})$.\\

\end{enumerate}

\item Each $\K-$algebra automorphism $\varphi$ of $B[w]$ is triangular. In particular, any $\K-$algebra automorphism $\varphi$ of $B[w]$ is given as follows:
\[
\varphi(x)=\varphi_{0}(x), \quad \varphi(y)=\varphi_{0}(y),\quad \varphi(s)=\varphi_{0}(s)
\]
and
\[
\varphi(t)=\varphi_{0}(t)\quad \varphi(w)=aw+b
\]
for some $\K-$algebra automorphism $\varphi_{0}$ of $B$, and $b$ in the center of $B$, and $a\in \K^{\ast}$ . \\

\item If $\Z\subseteq \K$, then ${\rm LNDer}(B)=\{0\}$.

\end{enumerate}
\end{thm}

{\bf Proof:} Recall that $xy-pyx=1$, and $st=qts$, and $p\neq 1$, and $q\neq 1$. It is easy to see that each affine $\K-$algebra automorphism $\varphi$ of $B$ has to be linear with respect to $x,y, s, t$. As a result, we have the following
\begin{eqnarray*}
\varphi(x)&=&a_{11} x+a_{12}y+a_{13} s+a_{14} t,\\
\varphi(y)&=&a_{21} x+a_{22} y+a_{23} s+a_{24} t,\\
\varphi(s)&=&a_{31} x+a_{32} y+a_{33} s+a_{34} t,\\
\varphi(t)&=&a_{41} x+a_{42} y+a_{43} s+a_{44} t
\end{eqnarray*}
where $a_{ij}\in \K$ for $i, j=1,2, 3, 4$.
Applying $\varphi$ to $xy-pyx=1$, we further have the following 
\begin{eqnarray*}
1&=& (a_{11}x+a_{12}y+a_{13}s+a_{14}t)(a_{21}x+a_{22}y+a_{23}s+a_{24}t)\\
&\quad & -p(a_{21}x+a_{22}y+a_{23}s+a_{24}t)(a_{11}x+a_{12}y+a_{13}s+a_{14}t)\\
&=& a_{11} a_{21}(1-p)x^{2}+(a_{11}a_{23}+a_{21}a_{13})(1-p)xs+(a_{11}a_{24}\\
&\quad &+a_{21}a_{14})(1-p)xt +a_{12}a_{22}(1-p)y^{2}+(a_{12}a_{23}+a_{22}a_{13})(1-p)ys\\
&\quad & +a_{13}a_{23}(1-p)s^{2}+(a_{14}a_{22}+a_{12}a_{24})(1-p)ty+a_{14}a_{24}(1-p)t^{2}\\
&\quad & +a_{11}a_{22}(xy-pyx)+a_{12}a_{21}(yx-pxy)\\
&\quad & +[a_{13}a_{24}(1-pq^{-1})+a_{14}a_{23}(q^{-1}-p)]st.
\end{eqnarray*}
Notice that $p\neq 1$. If $p=-1$, then the above equation implies that either
\[
a_{11}=a_{13}=a_{14}=0,\, a_{21}=a_{12}^{-1}\neq 0,\,
a_{21}=a_{12}^{-1}, \, a_{22}=a_{23}=a_{24}=0;
\]
or
\[
a_{11}=a_{22}^{-1}\neq 0,\, a_{12}=a_{13}=a_{14}=0,\, 
a_{21}=a_{23}=a_{24}=0, \,a_{22}=a_{11}^{-1}.
\]
If $p\neq -1$, then above equation implies that 
\[
a_{11}=a_{22}^{-1}\neq 0, \, a_{12}=a_{13}=a_{14}=0,\,
a_{21}=a_{23}=a_{24}=0, \,a_{22}=a_{11}^{-1}.
\]
As a result, if $p\neq -1$, then we have the following: 
\[
\varphi(x)=\alpha x, \quad \varphi(y)=\alpha^{-1} y
\]
for some $\alpha \in \K^{\ast}$. If $p=-1$, then we have that either 
\[
\varphi(x)=\alpha x,\quad \varphi(y)=\alpha^{-1}y
\]
for some $\alpha \in \K^{\ast}$ or
\[
\varphi(x)=\alpha y, \quad \varphi(y)=\alpha^{-1}x
\]
for some $\alpha \in \K^{\ast}$. 

In either case, we shall have that $a_{31}=a_{32}=a_{41}=a_{42}=0$ due to the fact that $xt=tx, xs=sx$ and $yt=ty, ys=sy$ and $p\neq 1$. Also note that $st=qts$ and $q\neq 1$. If $q=-1$, then we shall have that either
\[
a_{33}\neq 0, \quad a_{34}=0, \quad a_{43}=0,\quad a_{44}\neq 0;
\]
or 
\[
a_{33}=0,\quad a_{34}\neq 0, \quad a_{43}\neq 0, \quad a_{44}=0.
\]
If $q\neq -1$, we shall have that 
\[
a_{33}\neq 0,\quad a_{34}=0, \quad a_{43}=0, \quad a_{44}\neq 0
\]
As a result, if $q=-1$, then we have that either
\[
\varphi(s)=\beta s,\quad \varphi(t)=\gamma t
\]
for some $\beta, \gamma \in \K^{\ast}$ or 
\[
\varphi(s)=\beta t, \quad \varphi(t)=\gamma s.
\]
for some $\beta, \gamma \in \K^{\ast}$. Otherwise, if $q\neq -1$, then we have that 
\[
\varphi(s)=\beta s, \quad \varphi(t)=\gamma t
\]
for some $\beta, \gamma \in \K^{\ast}$. Therefore, Part (1) of the statement follows.

Parts (2) and (3) follow from the the fact the algebra $B$ is in the category ${\bf Af}$ as defined in \cite{CPWZ1}, when both $p$ and $q$ are roots of unity. And the method of modulo p can be used to reduce the problem to the situation where both $p$ and $q$ are roots of unity, and thus the discriminant of $B$ is dominating. Now the rest of the results follow from {\bf Theorem 1.13} in \cite{CPWZ2}. We will not repeat the details here.
\qed

\vskip .5in

\noindent
{\bf Acknowledgements:} The author would like to thank James Zhang for some useful communications regarding the discriminant method and the cancellation property.

\end{document}